\documentclass[11pt, a4paper, reqno]{amsart}
\usepackage{amssymb, array, amsmath, amscd, pdfpages, enumerate, amsthm, fixltx2e, setspace, hyperref,mathtools}

\usepackage[utf8]{inputenc}
\usepackage{amssymb}
\usepackage{bm}
\usepackage{graphicx}

\DeclareMathOperator{\spsup}{sp \;sup}

\newcommand*{\R}{{\mathbb{R}}}     
\newcommand*{\Z}{{\mathbb{Z}}}     
     
\newcommand*{\N}{{\mathbb{N}}}

\newcommand*{\Abs}[2][default]{\ifthenelse{\equal{#1}{default}}{\left\lvert#2\right\rvert}{\ldelim{#1}{\lvert}#2\rdelim{#1}{\rvert}}}
\newcommand*{\Norm}[2][default]{\ifthenelse{\equal{#1}{default}}{\left\lVert#2\right\rVert}{\ldelim{#1}{\lVert}#2\rdelim{#1}{\rVert}}}

\newcommand*{\Iprod}[3][default]{\ifthenelse{\equal{#1}{default}}{\left\langle#2,#3\right\rangle}{\ldelim{#1}{\langle}#2,#3\rdelim{#1}{\rangle}}}
\newcommand*{\Dualpair}[3][default]{\ifthenelse{\equal{#1}{default}}{\left\langle#2,#3\right\rangle}{\ldelim{#1}{\langle}#2,#3\rdelim{#1}{\rangle}}}

\newcommand*{\ddb}[2][1]{\ifthenelse{\equal{#1}{1}}{\frac{d}{d#2}}{\frac{d^{#1}}{d#2^{#1}}}}
\newcommand*{\pd}[3][1]{\ifthenelse{\equal{#1}{1}}{\frac{\partial{#2}}{\partial{#3}}}{\frac{\partial^{#1}{#2}}{\partial#3^{#1}}}}

\newcommand*\lenv{{\hbox{\raisebox{-.15ex}{\rotatebox[origin=c]{50}{$\smallsmile$}}\kern-8.65pt\rotatebox[origin=c]{-25}{$\smallsetminus$}}}}
\newcommand*\uenv{{\hbox{\raisebox{-.0ex}{\rotatebox[origin=c]{-45}{$\smallfrown$}}\kern-5.6pt\raisebox{.2ex}{\rotatebox[origin=c]{-5}{\scriptsize\slash}}}}\,\kern+1.5pt}

\newcommand*\llenv{{\hbox{\raisebox{-.15ex}{\rotatebox[origin=c]{50}{$\smallsmile$}}\kern-7.15pt\rotatebox[origin=c]{-25}{$\smallsetminus$}}}}
\newcommand*\uuenv{{\hbox{\raisebox{-.0ex}{\rotatebox[origin=c]{-45}{$\smallfrown$}}\kern-4.6pt\raisebox{.2ex}{\rotatebox[origin=c]{-5}{\scriptsize\slash}}}}\,\kern+1.5pt}

\newcommand{\sleq}{\preccurlyeq}
\newcommand{\sgeq}{\succcurlyeq}
\newcommand*{\Gg}{{\mathcal{G}}}

\newcommand*{\Uu}{{\mathcal{U}}} 

\newcommand*{\Ww}{{\mathcal{W}}}

\newcommand*{\lupp}[1]{\kern0pt^{\kern0pt u \kern0pt}\kern1pt#1}
\newcommand*{\llow}[1]{\kern0pt^{\kern0pt l \kern0pt}\kern1pt#1}
\newcommand*{\rupp}[1]{#1\kern1pt^{\kern0pt u \kern0pt}\kern0pt}
\newcommand*{\rlow}[1]{#1\kern1pt^{\kern0pt l \kern0pt}\kern0pt}
\newcommand*{\ul}[1]{\kern0pt^{\kern0pt u \kern0pt}\kern0pt#1\kern1pt^{\kern0pt l \kern0pt}\kern0pt}
\newcommand*{\lu}[1]{\kern0pt^{\kern0pt l \kern0pt}\kern0pt#1\kern0pt^{\kern0pt u \kern0pt}\kern0pt}

\DeclareMathOperator{\mmin}{Min}
\DeclareMathOperator{\mmax}{Max}
\newcommand{\mini}[1]{\mmin [\kern1pt #1 \kern1pt ]}
\newcommand{\maxi}[1]{\mmax [\kern1pt #1 \kern1pt ]}

\newcommand*{\bproofname}{Proof}
\newenvironment{bproof}[1][\bproofname]{\begin{proof}[#1]}{\end{proof}}

\newtheorem{thm}{Theorem}[section]

\newtheorem{lemma}[thm]{Lemma}
\newtheorem{cor}[thm]{Corollary}
\theoremstyle{definition}
\newtheorem{defn}[thm]{Definition}

\newtheorem{remark}[thm]{Remark}
\newtheorem{example}[thm]{Example}

\numberwithin{equation}{section}

\begin{document}

\title[Partially ordered semigroups and groups with two orderings]{Partially ordered semigroups and groups with two mixed partial orderings}

\thispagestyle{plain}

\author{Jani Jokela}
\address[J. Jokela]{}
\email{jani.jokela@tuni.fi}
%\email{jani.m.jokela@gmail.com}
%
% Organisation:
% Mathematics Research Centre, (Faculty of Information Technology and Communication Sciences,) 
% Tampere University, Tampere, Finland
%
%\address[J. Jokela]{Department of Mathematics, Tampere University of Technology, P.O.\ Box 553, 33101 Tampere, Finland}
%

\begin{abstract}
A mixed lattice is a partially ordered set with two mixed partial orderings that are linked by asymmetric upper and lower envelopes. These notions generalize the join and meet operations of a lattice.  
In the present paper, we study different types of partially ordered semigroups with two mixed orderings, and investigate their relationship to subsemigroups of mixed lattice groups, which are partially ordered groups with a similar order structure.  
We also consider Archimedean orderings, and we show that  
elements of finite order cannot exist in a rather general class of Archimedean mixed lattice groups. Moreover, we give an example of a non-Archimedean mixed lattice group that contains an element of finite order.  
\end{abstract}


\maketitle

\section{Introduction}

A \emph{mixed lattice} is a certain type of generalization of a lattice. Mixed lattices are partially ordered sets with two partial orderings that are linked together by asymmetric upper and lower envelopes formed with respect to the two partial orderings. These mixed envelopes are generalizations of the join and meet operations in lattice theory. Equivalently, a mixed lattice can be viewed as a non-commutative and non-associative algebra. These rather general ideas were introduced in \cite{jj5}. 
However, the origins of the mixed lattice theory dates back several decades to the work of Arsove and Leutwiler \cite{ars2, ars3, ars4, ars}. They developed an axiomatization of potential theory based on partially ordered semigroups with two mixed partial orderings. The study of mixed lattice theory was later carried over to groups by Eriksson-Bique \cite{eri1, eri}, and more recently, to a vector space setting \cite{jj1, jj2, jj3}. Nevertheless, the mixed lattice theory has not yet been very widely studied, especially in the group setting, and current knowledge on these structures remains somewhat limited. Indeed, even some of the more elementary aspects of the theory are not very well understood. The purpose of this paper is therefore to fill some of the gaps in the basic theory of ordered commutative semigroups and groups with two mixed partial orders.

In the next section, we first study different types of 
semigroups with two mixed orderings. We collect the main facts from previous studies and supplement them with new results. In Section \ref{sec:groups}, we study the relationships between mixed lattice groups and the semigroup notions of Section \ref{sec:semig}. Mixed lattice groups can be divided into so-called regularity classes based on the generality of the assumptions regarding the two partial orderings. We review these classes and their basic properties and characterizations in Section \ref{sec:groups}. 
In Section \ref{sec:archi}, we consider mixed lattice groups with the Archimedean property.  
In particular, we study the existence of elements of finite order. It is well known that in lattice ordered groups every element has infinite order. In contrast, a mixed lattice group may contain elements of finite order  
if at least one of the partial orderings is not an Archimedean order. We will show that there is a rather general class of Archimedean mixed lattice groups in which elements of finite order cannot exist.

%----------------------

\section{Semigroups with two mixed partial orderings}\label{sec:semig}

We recall that a commutative semigroup $(S,+)$ together with a partial ordering $\leq$ is called a 
\emph{partially ordered semigroup}, if  
\begin{equation}\label{pogroup}
u\leq v \; \implies \; u+w\leq  v+w  
\end{equation}
holds for all $u,v,w\in S$.

A partially ordered commutative semigroup $(S,+,\leq)$ is called \emph{naturally ordered} if 
\begin{equation}\label{natural}
u\leq v \; \iff \; v=u+w   \quad \text{ for some } \, \,w\in S.
\end{equation}

If $S$ is naturally ordered and $S$ contains a neutral element (denoted by $0$) then $S$ is 
\emph{positively ordered}, that is, $x\geq 0$ holds for all $x\in S$.

Next we shall assume that $S$ is a 
commutative partially ordered semigroup with a neutral element $0$  
and two partial orderings $\leq$ and $\sleq$ (i.e. we assume that both orderings satisfy \eqref{pogroup}). In this context, the partial order $\leq$ is traditionally called the \emph{initial order} and $\sleq$ is called the \emph{specific order}.

For every $x,y\in S$ the \emph{mixed upper envelope} is defined by 
\begin{equation}\label{upperenv}
u\uenv v\,=\,\min \,\{\,w\in S: \; w\sgeq u \; \textrm{ and } \; w\geq v \,\},
\end{equation}
and the \emph{mixed lower envelope} is defined by  
\begin{equation}\label{lowerenv}
u\lenv v \,=\,\max \,\{\,w\in S: \; w\sleq u \; \textrm{ and } \; w\leq v \,\}, 
\end{equation}
where the minimum and maximum (if they exist) are taken with respect to the initial order $\leq$. 
Here, and everywhere in this paper, 
it is understood that the minimum and maximum of a set are always taken with respect to the partial order $\leq$.

It follows immediately from the definitions that 
\begin{equation}\label{p0}
x\lenv y \sleq x\sleq x\uenv y \; \textrm{ and } \;  x\lenv y \leq y\leq x\uenv y
\end{equation}
and
\begin{equation}\label{p4}
x\sleq u \; \textrm{ and } \; y\leq v \; \implies \; x\uenv y \leq u\uenv v \; \textrm{ and } \; x
\lenv y \leq  u\lenv v 
\end{equation}
hold for all $x,y,u,v\in S$ assuming all the mixed envelopes exist. Furthermore, \eqref{p0} and \eqref{p4} imply that 
\begin{equation}\label{p5b}
x\sleq y \, \iff \, x\uenv y = y \, \iff \, y\lenv x = x
\end{equation}
holds in $S$. Proofs of these facts can be found in \cite{jj5}.

We now have the following definition.

\begin{defn}\label{weak_mlg}
Let $\Uu=(S,+,\leq,\sleq)$ be a partially ordered commutative 
semigroup with two partial orderings and a neutral element $0$. We assume that $\Uu$ is naturally ordered with respect to $\sleq$, and that $\leq$ is any partial ordering that satisfies the strong cancellation property  
\begin{equation}\label{canc1}
u\leq v \quad \iff \quad u+w \leq v+w \quad \textrm{for all } u,v,w\in \Uu.
\end{equation}    
Then 
\begin{enumerate}[(i)]
\item
$\Uu$ is called an \emph{$A$-structure} if the element $x\uenv y$ exists in $\Uu$ for all $x,y\in \Uu$.
\item
$\Uu$ is called a \emph{positive A-structure} if $\Uu$ is an $A$-structure such that $x\geq 0$ for all $x\in \Uu$.
\item
$\Uu$ is called a \emph{weak mixed lattice semigroup} if the elements $x\uenv y$ and $x\lenv y$ exist in $\Uu$ for all $x,y\in \Uu$.
\item
$\Uu$ is called a \emph{mixed lattice semigroup} if $\Uu$ is a weak mixed lattice semigroup such that $x\uenv y + y\lenv x = x+y$ for all $x,y\in \Uu$.
\end{enumerate}
\end{defn}

We note that, with the above assumptions, the strong cancellation property \eqref{canc1} for $\leq$ implies that the partial order $\sleq$ has the same property. Indeed, $x+z\sleq y+z$ is equivalent to $y+z=x+z+w$ for some $w\in\Uu$. By \eqref{canc1} this is equivalent to $y=x+w$, or $x\sleq y$. 

The notion of $A$-structure was introduced in \cite{ars} and  
it was studied further in \cite{eri}. We should remark that 
the authors in \cite{ars} gave a definition of an $A$-structure that was slightly different than the definition used in \cite{eri} and in the present paper. In \cite{ars}, %the semigroup $\Uu$ 
an $A$-structure was always assumed to be positive but it was not assumed to satisfy the cancellation property \eqref{canc1}. Moreover, the results in \cite{eri} were presented for positive $A$-structures, but many of the proofs actually work as such without the positivity assumption. In this paper we assume the cancellation property but positivity is assumed only when necessary.

The following characterizations are from \cite{eri} and \cite{ars}.

\begin{thm}\label{karakt11}
Let $\Uu=(S,+,\leq,\sleq)$ be a partially ordered commutative 
semigroup with two partial orderings and the neutral element such that $\leq$ satisfies \eqref{canc1} and $\Uu$ is naturally ordered with respect to $\sleq$. Then 
\begin{enumerate}[(a)]
\item
$\Uu$ is an $A$-structure if and only if $m=\min \{w\in \Uu: u\leq w+v\}$ exists for all $u,v\in \Uu$. The minimum element satisfies $m+v=v\uenv u$.
\item
If $\Uu$ is a positive $A$-structure then $\Uu$ is a weak mixed lattice semigroup if and only if 
$m=\min \{w\in \Uu: u\leq w+v, \; w\sleq u\}$ exists for all $u,v\in \Uu$. In this case, the minimum element satisfies $m+ u\lenv v=u$.  
\item
$\Uu$ is a mixed lattice semigroup if and only if $m=\min \{w\in \Uu: u\leq w+v\}$ exists for all $u,v\in \Uu$ and it satisfies $m\sleq u$. In this case, $m$ satisfies $m+v=v\uenv u$ and $m+ u\lenv v=u$. 
\end{enumerate}
\end{thm}

\begin{remark}
In \cite{eri}, the existence of the minimum element in Theorem \ref{karakt11}(a) was taken as the definition of an $A$-structure. In the present paper, we have chosen to use the equivalent statement of Definition \ref{weak_mlg} instead as the definition because it resembles the definition of an upper semilattice, and therefore appears to be a simpler and more illuminating than the condition of Theorem \ref{karakt11}(a).
\end{remark}

It was proved in \cite[Corollary 2.3]{eri} that in every $A$-structure the mixed upper envelope is translation invariant, that is, if $\Uu$ is an $A$-structure then 
\begin{equation}\label{astruct-eq1}
(u+a)\uenv (v+a)=u \uenv v +a
\end{equation}
holds for all $u,v,a\in\Uu$.

It was also proved in \cite[Theorem 2.5]{eri} that if $\Uu$ is a weak mixed lattice semigroup then the translation invariance of the mixed lower envelope is equivalent to $\Uu$ being a mixed lattice semigroup.

Moreover, an $A$-structure $\Uu$ is a mixed lattice semigroup if and only if the inequality $u\uenv v \sleq u+v$ holds for all $u,v\in \Uu$ (\cite[Theorem 1.1]{ars}). We now give a similar characterization for positive $A$-structures. 

\begin{thm}\label{positive_astr_char}
An $A$-structure $\Uu$ is a positive $A$-structure if and only if the inequality $v\uenv u \leq u+v$ holds for all $u,v\in \Uu$.
\end{thm}

\begin{bproof}
Let $\Uu$ be a positive $A$-structure with $u,v\in \Uu$ and define $V=\{w\in \Uu: u\leq w+v\}$. Then by Theorem \ref{karakt11} the set $V$ has the smallest element $m$ that satisfies $m+v=v\uenv u$. Since $v\geq 0$, we have $u\leq u+v$, and so $u\in V$. Hence $m\leq u$ and $v\uenv u = m+v \leq u+v$.  

Conversely, if the given inequality holds in $\Uu$ then setting $u=0$ gives $v=v+0\geq v\uenv 0 \geq 0$ 
for all $v\in \Uu$.
\end{bproof}

Next we turn to weak mixed lattice semigroups. 
First, it follows from Theorem \ref{karakt11} that if $\Uu$ is a weak mixed lattice semigroup then 
\begin{equation}\label{wmlsg-ineq_1}
v\uenv u + u\lenv v \leq u+v
\end{equation}
holds for all $u,v\in\Uu$. 
Indeed, if we put $V=\{w\in \Uu: u\leq w+v\}$ and $W=\{w\in \Uu: u\leq w+v, \, w\sleq u\}$ then obviously $W\subseteq V$. Hence $m_1=\min V \leq \min W=m_2$, and we can apply the equalities given in Theorem \ref{karakt11}(a) and (b) to obtain
$$
v\uenv u + u\lenv v = m_1+v + u\lenv v \leq m_2 +v + u\lenv v = u+v. 
$$ 
If equality holds in \eqref{wmlsg-ineq_1} for all $u,v\in \Uu$ then $\Uu$ is a mixed lattice semigroup.

The next theorem extends the result  
given in Theorem \ref{karakt11}(b). It also shows that, in a certain sense, 
the inequality \eqref{wmlsg-ineq_1} characterizes weak mixed lattice semigroups.

\begin{thm}\label{wmlsg_karakt2}
If $\Uu$ is a positive $A$-structure then the following conditions are equivalent.
\begin{enumerate}[(a)]
\item
$\Uu$ is a weak mixed lattice semigroup. 
\item
The element $m=\min \{w\in \Uu: u\leq w+v, \; w\sleq u\}$ exists for all $u,v\in \Uu$.
\item
The element
$$
M=\max \{w\in \Uu: v\uenv u + w\leq u+v, \; w\sleq u\} 
$$
exists for all $u,v\in \Uu$. In fact, $M=u\lenv v$.
\end{enumerate}
\end{thm}

\begin{bproof}
Let us denote the sets by 
$V=\{w\in \Uu: u\leq w+v, \; w\sleq u\}$ and $W=\{w\in \Uu: v\uenv u + w\leq u+v, \; w\sleq u\}$. 
First we note that the positivity of $\Uu$ guarantees that the sets $V$ and $W$ are nonempty. Indeed, if $u,v\in \Uu$ then $u\sleq u$ and $0\leq v$, so $u\leq u+v$. This shows that $u\in V$. On the other hand, Theorem \ref{positive_astr_char} implies that $0\in W$. 

The equivalence of (a) and (b) was established in \cite[Theorem 2.4]{eri}. We will prove the equivalence (a)$\iff$(c). 
If (a) holds then $u\lenv v$ exists in $\Uu$ and $u\lenv v \in W$, by \eqref{wmlsg-ineq_1}. Let $w\in W$. Then $w\sleq u$ and 
\begin{equation}\label{apuineq1}
v\uenv u + w \leq u+v \leq v\uenv u +v,
\end{equation}
and it follows by the cancellation property that $w\leq v$. Hence $w\leq u\lenv v$, and this shows that $u\lenv v = \max W$.

Next, assume that (c) holds. If $M=\max W$ then obviously $M\sleq u$. Also, the inequalities \eqref{apuineq1} hold for $w=M$, and so $M\leq v$. We will show that $M=u \lenv v$. For this, choose some $w\in \Uu$ such that $w\leq v$ and $w\sleq u$. The latter inequality implies that $u=w+ m$ for some $m\in \Uu$. Then $v\sleq v+ m$ and $u=w+m\leq v+m$, and consequently, $v\uenv u \leq v+m$, by \eqref{p4}. Hence we have
$v\uenv u + w \leq v+ m+w = v+u$. This shows that $w\in W$, and so $w\leq M$. Hence $M=u\lenv v$, proving that $\Uu$ is a weak mixed lattice semigroup.
\end{bproof}

Theorem \ref{wmlsg_karakt2} gives conditions for a positive $A$-structure to be a weak mixed lattice semigroup. On the other hand, every weak mixed lattice semigroup is a positive $A$-structure.  
This is a rather immediate consequence of the definitions.

\begin{thm}\label{wmlsg_is_positive}
Every weak mixed lattice semigroup is a positive $A$-structure.
\end{thm}

\begin{bproof}
If $\Uu$ is a weak mixed lattice semigroup and $u\in \Uu$ then  
the element $0\lenv u$ exists in $\Uu$, and so $0\sleq 0\lenv u \sleq 0$, which implies that $0\lenv u=0$. Hence, $u\geq 0\lenv u =0$ and $\Uu$ is a positive $A$-structure.
\end{bproof}

\begin{example}\label{a_structure}
Let $S=\{0,1,2,3,4,\ldots\}$, let $\sleq$ be the usual ordering of integers and let $\leq$ be the ordering corresponding to the positive cone $S_{p}=\{3k:k\in V\}$. Then $\Uu=(S,+,\leq,\sleq)$ is an $A$-structure 
since the element $x\uenv y$ exists in $\Uu$ for all $x,y\in \Uu$, but not a positive $A$-structure, and hence not a weak mixed lattice semigroup. Note that now $z=y\lenv x$ is the largest number such that $y\sleq z$ and $z$ is of the form $z=x-3k$, where $k\in V$. Thus, for example, the element $1\lenv 2$ does not exist in $\Uu$.
\end{example}

An example of a weak mixed lattice semigroup that is not a mixed lattice semigroup was given in 
\cite[Example 2.3]{jj1}. We shall further discuss the same setting in Example \ref{weak_mlsg_ex2}.

We have already mentioned a couple of conditions under which an $A$-structure is a mixed lattice semigroup. For reference, we collect a few of them in the next theorem. 
The following characterizations were proved in \cite[Theorem 2.5]{eri}. The equivalence of (a) and (b) and the implication (a)$\implies$(c) of the following theorem were also proved earlier in \cite{ars} with different proofs than those given in \cite{eri}. We give yet another proof below based on our results thus far.

\begin{thm}\label{mlsgchar1}
Let $\Uu$ be an $A$-structure. Then the following statements are equivalent.
\begin{enumerate}[(a)]
\item
$\Uu$ is a mixed lattice semigroup. 
\item 
The inequality $v\uenv u\sleq u+v$ holds for all $u,v\in\Uu$. 
\item
The mixed lower envelopes exist in $\Uu$ and $(u+a)\lenv (v+a)=u \lenv v +a$ holds for all $u,v,a\in\Uu$. 
\end{enumerate}
\end{thm}

\begin{bproof}
If (a) holds then $u\lenv v$ exists and $u\lenv v\sgeq 0$ for all $u,v\in \Uu$. Hence $v\uenv u \sleq v\uenv u + u\lenv v = u+v$, where the last equality holds by the definition of a mixed lattice semigroup.

Suppose that (b) holds. First we note that $v\sleq v\uenv 0 \sleq v+0=v$, and so $0\leq v\uenv 0=v$ for all $v\in \Uu$. Hence $\Uu$ is a positive $A$-structure. 
Next, let $W=\{w\in \Uu: v\uenv u + w\leq u+v, \; w\sleq u\}$. Since $v\uenv u \sleq u+v$ there exists an element $m\in\Uu$ such that $v\uenv u +m= u+v$. Moreover, from $v\sleq v\uenv u$ it follows that 
$v+m\sleq v\uenv u + m=u+v$. Canceling $v$ gives $m\sleq u$. Hence $m\in W$ and $m$ is clearly the largest element in $W$ since it gives equality in $v\uenv u + m\leq u+v$. Thus $m=u\lenv v$ by Theorem \ref{wmlsg_karakt2}(c),  
and so we actually have $v\uenv u + u\lenv v = u+v$, and hence $\Uu$ is a mixed lattice semigroup. This establishes the equivalence of (a) and (b).

Assuming (c) holds,  
let $\Uu$ be a weak mixed lattice semigroup such that $(u+a)\lenv (v+a)=u \lenv v +a$ holds for all $u,v,a\in\Uu$. Noting that $u+v\uenv u\sgeq u+v$ and $v+v\uenv u \geq u+v$, then using \eqref{astruct-eq1} and \eqref{p4} we obtain 
$$
v\uenv u + u\lenv v = (u+v\uenv u)\lenv (v+v\uenv u) \geq (u+v)\lenv (u+v) =u+v.
$$
The reverse inequality holds by \eqref{wmlsg-ineq_1}. Hence $v\uenv u + u\lenv v = u+v$, and $\Uu$ is a mixed lattice semigroup. Thus (c) implies (a).  

Conversely, if $\Uu$ is a mixed lattice semigroup then 
$$
(a+v)\uenv(a+u)+(a+u)\lenv(a+v) = u+v+a+a
$$
holds for all $u,v,a\in\Uu$. Similarly, $u+v=v\uenv u + u\lenv v$. Moreover, by \eqref{astruct-eq1} we have $(a+v)\uenv (a+u)=v \uenv u +a$. Substituting these into the above equation gives
$$
v \uenv u +a +(a+u)\lenv(a+v) = v\uenv u + u\lenv v + a+a.
$$
Canceling the terms $v\uenv u +a$ yields the desired result and completes the proof. (The proof of this last implication is due to Eriksson-Bique \cite{eri}.) 
\end{bproof}

Next we introduce a kind of a dual notion of an $A$-structure with many similar properties.

\begin{defn}\label{b-structure}
Let $\Ww=(S,+,\leq,\sleq)$ be a partially ordered Abelian semigroup with two partial orderings $\leq$ and $\sleq$, and a neutral element $0$. We assume that $\Ww$ is naturally ordered with respect to $\leq$, and that $\sleq$ is any partial ordering that satisfies the strong cancellation property  
\begin{equation}\label{canc2}
u\sleq v \quad \iff \quad u+w \sleq v+w \quad \textrm{for all } u,v,w\in \Ww.
\end{equation}
If the mixed upper envelope $v\uenv u$ exists in $\Ww$ for all $u,v\in \Ww$ then $\Ww$ is called a \emph{$B$-structure}. 
If  $\Ww$ is a $B$-structure and $u\sgeq 0$ for all $u\in \Ww$ then $\Ww$ is called a \emph{positive $B$-structure}.
\end{defn}

As with $A$-structures, it again follows from \eqref{canc2} that 
the strong cancellation property holds also for 
the partial order $\leq$. It is also clear from the definitions that if $\Ww$ is both an $A$-structure and a $B$-structure, then the two partial orderings must coincide, and $\Ww$ is then an upper semilattice.

At first glance, it may appear as if the definition of a $B$-structure is the same as the definition of an $A$-structure with merely the roles of the two orderings interchanged, but this is not quite true since the minimum in the definition of the mixed upper envelope is still taken with respect to the partial order $\leq$.

The proofs of the following results are similar to the proofs of the corresponding results for $A$-structures.  
Nevertheless, we will give full details for the next theorem  
to indicate the minor differences in the proofs.

\begin{thm}\label{karakt22}
Let $\Ww=(S,+,\leq,\sleq)$ be a partially ordered commutative 
semigroup with two partial orderings and the neutral element such that $\sleq$ satisfies \eqref{canc2} and $\Ww$ is naturally ordered with respect to $\leq$.  If $u,v\in \Ww$ then the following statements are equivalent.
\begin{enumerate}[(a)]
\item
The mixed upper envelope $u\uenv v$ exists in $\Ww$.
\item
The element $m=\min \{w\in \Ww: u\sleq w+v\}$ exists in $\Ww$. 
\end{enumerate}
Moreover, the minimum element in (b) satisfies $m+v=u\uenv v$. In particular, $\Ww$ is a $B$-structure if and only if $m=\min \{w\in \Ww: u\sleq w+v\}$ exists for all $u,v\in \Ww$.
\end{thm}

\begin{bproof}
Let $V=\{w\in \Ww: u\sleq w+v\}$. If $u\uenv v$ exists then $v\leq u\uenv v$, and since $\Ww$ is naturally ordered with respect to $\leq$, there exists some $m \in \Ww$ such that $m+v= u\uenv v \sgeq u$. Hence $m \in V$. 
Let $w\in V$. Then $u\sleq w+v$ and $v\leq w+v$ (since $w\geq 0$). Consequently, $w+v\geq u\uenv v = m+v$, and it follows by the cancellation property that $w\geq m$. This shows that $m=\min V$.

Conversely,  assume that $m=\min V$ exists. Then $u\sleq m+v$ and $v\leq m+v$, since $m\geq 0$. We will show that $m+v=u\uenv v$. Choose any $w\in \Ww$ such that $u\sleq w$ and $v\leq w$. Then, since $\Ww$ is naturally ordered with respect to $\leq$, there exists some $v' \in \Ww$ such that $v'+v=w\sgeq u$. Hence $v' \in V$ and it follows that $m\leq v'$. Thus we have $m+v\leq v'+v=w$, and this shows that $m+v=u\uenv v$.
\end{bproof}

It follows from the preceding theorem that the translation invariance of the mixed upper envelopes holds also in every $B$-structure. That is, if $\Ww$ is a $B$-structure then 
\begin{equation}\label{bstruct-eq1}
(u+a)\uenv (v+a)=u \uenv v +a
\end{equation}
holds for all $u,v,a\in\Ww$. This can be proved as follows. By the cancellation property \eqref{canc2}, the inequalities $u\sleq w+v$ and $u+a\sleq w+v+a$ are equivalent, and hence the sets 
$$
\{w\in \Ww: u\sleq w+v\} \quad \text{and} \quad \{w\in \Ww: u+a\sleq w+v+a\}
$$
are equal. Thus, by Theorem \ref{karakt22}, they both have the same minimum element $m$ that satisfies the equations
$$
m+v=u\uenv v \quad \text{and} \quad m+v+a=(u+a)\uenv (v+a). 
$$
Adding $a$ to the first equation yields \eqref{bstruct-eq1}.

The next two results are counterparts of Theorems \ref{positive_astr_char} and \ref{wmlsg_karakt2}, respectively. We will omit the proofs since they are similar to the proofs of the aforementioned theorems.

\begin{thm}\label{positive_bstr_char}
A $B$-structure $\Ww$ is a  %$(\sleq)$-
positive $B$-structure if and only if the inequality $u\uenv v \leq u+v$ holds for all $u,v\in \Ww$. 
\end{thm}

%\begin{bproof}
%Let $\Ww$ be a positive $B$-structure with $u,v\in \Ww$ and define $V=\{w\in \Ww: u\sleq w+v\}$. Then by Theorem %\ref{karakt22} the set $V$ has the smallest element $m$ that satisfies $m+v=u\uenv v$. Since $v\sgeq 0$, we have %$u\sleq u+v$, and so $u\in V$. Hence $u\uenv v = m+v \leq u+v$. %$v\uenv u -v\leq u$ and the result follows. 

%Conversely, if the given inequality holds in $\Ww$ then setting $u=0$ gives $v=v+0\sgeq 0\uenv v \sgeq 0$
%, and so $v\geq 0$ 
%for all $v\in \Ww$.
%\end{bproof}

\begin{thm}\label{karakt23}
Let $\Ww=(S,+,\leq,\sleq)$ be a positive $B$-structure.   
If $u,v\in \Ww$ then the following statements are equivalent.
\begin{enumerate}[(a)]
\item
The mixed lower envelope $v\lenv u$ exists in $\Ww$.
\item
The element $m=\min \{w\in \Ww: u\sleq w+v, \; w\leq u\}$ exists.  
\item
The element
$M=\max \{w\in \Ww: u\uenv v + w\sleq u+v, \; w\leq u\}$ exists. 
\end{enumerate}
Moreover, the minimum element in (b) satisfies $m+ v\lenv u=u$, and the maximum element in (c) is $M=v\lenv u$. 
\end{thm}

The preceding theorem implies that if $\Ww$ is a $B$-structure such that the mixed lower envelope $u\lenv v$ exists for all $u,v\in \Ww$ then the inequality 
\begin{equation}\label{wmlsg-ineq_222}
v\uenv u + u\lenv v \sleq u+v
\end{equation}
holds in $\Ww$.

Many concrete examples of $A$-structures and mixed lattice semigroups are constructed as follows. Let $\Ww$ be a positively partially ordered semigroup that satisfies the strong cancellation law \eqref{canc1}, and let $\Uu$ be a subsemigroup of $\Ww$ containing the neutral element. Let $\sleq$ be the partial order whose positive set is the subsemigroup $\Uu$, and let $\leq$ be the partial order of $\Ww$ restricted to $\Uu$. Then $(\Uu,\leq,\sleq)$ is a partially ordered semigroup with two partial orderings such that $\Uu$ has the cancellation property \eqref{canc1} and $\Uu$ is naturally ordered with respect to $\sleq$. If the mixed upper envelope $x\uenv y$ exists in $\Uu$ then $(\Uu,\leq,\sleq)$ is a positive $A$-structure. 

If a positive $A$-structure $\Uu$ arises from an oversemigroup $\Ww$ in this way, then the partial order $\sleq$ determined by $\Uu$ and, in certain cases, also the mixed upper envelopes of $\Uu$ can be extended to $\Ww$. If such extensions are possible then $\Ww$ is a $B$-structure.  

\begin{thm}
Let $(\Uu,\leq,\sleq)$ be a positive $A$-structure and assume that there exists an oversemigroup $\Ww$ of $\Uu$ such that $\Ww$ is naturally ordered with respect to $\leq$ and \eqref{canc1} holds in $\Ww$. Then the partial order $\sleq$ can be extended to $\Ww$ by defining 
\begin{equation}\label{sleq_ext}
a\sleq b \quad \text{in} \; \Ww \quad \text{if and only if} \quad b=a+w \quad \text{for some} \; w\in\Uu. 
\end{equation} 
If, in addition, for every $z\in \Ww$ there exists an element $w\in \Uu$ such that $z\sleq w$  
then the mixed upper envelopes of $\Uu$ can be extended to $\Ww$, and then $\Ww$ is a $B$-structure.  
\end{thm}

\begin{bproof}
We need to verify that the relation given by \eqref{sleq_ext} is indeed a partial ordering that satisfies \eqref{canc2}. Let $a,b,c\in \Ww$. 
Since $a=a+0$, it follows that $a\sleq a$ for all $a\in \Ww$. If $a\sleq b$ and $b\sleq a$ then $a=b+w$ and $b=a+v$ for some $v,w\in \Uu$. Then $a=a+w+v$, and by \eqref{canc1} we have $w+v=0$. Hence $w=0$ and $v=0$ (see \cite{hebi}, Corollary 1.11), and so $a=b$. Finally, if $a\sleq b$ and $b\sleq c$ then $b=a+w$ and $c=b+v$ for some $v,w\in \Uu$. Consequently, $c=a+w+v$ with $w+v\in\Uu$, so that $a\sleq c$. This shows that $\sleq$ is a partial ordering on $\Ww$.  Moreover, $\sleq$ satisfies the cancellation property \eqref{canc2} (see the remarks that follow Definition \ref{weak_mlg}).

Assume then that for every $z\in \Ww$ there exists an element $w\in \Uu$ such that $z\sleq w$ . 
If $a,b\in\Ww$ then $a\sleq x$ and $b\sleq y$ for some $x,y\in \Uu$. Hence $x=u+a$ and $y=v+b$ for some $u,v\in\Uu$, by \eqref{sleq_ext}. 
Since $\Uu$ is a positive $A$-structure, the element $(x+v)\uenv (y+u)$ exists in $\Uu$. Moreover, $b\geq 0$ and therefore we have 
$$
(x+v)\uenv (y+u)\geq y+u=u+v+b\geq u+v.
$$
Since $\Ww$ is naturally ordered with respect to $\leq$, there exists an element $m\in \Ww$ such that $u+v+m=(x+v)\uenv (y+u)$. To finish the proof we will show that this element $m$ equals $a\uenv b$. First we note that 
$$
u+v+m=(x+v)\uenv (y+u)\sgeq x+v=u+v+a,
$$
and by cancelling $u+v$ we obtain $a\sleq m$. Similarly, from 
$$
u+v+m=(x+v)\uenv (y+u)\geq y+u=u+v+b
$$
we get $b\leq m$. 
Assume then that $a\sleq w$ and $b\leq w$. Then $u+v+w\sgeq u+v+a=x+v$ and $u+v+w\geq u+v+b=y+u$, and hence by \eqref{p4} we have  
$$
u+v+w\geq (x+v)\uenv (y+u)=u+v+m.
$$
Consequently, $w\geq m$. This shows that $m=a\uenv b$. 
\end{bproof}

\section{Mixed lattice groups}\label{sec:groups}

In this section we study the relationship between the results presented in the preceding section and the structure of mixed lattice groups. Let $G=(\Gg, +,\leq, \sleq)$ be a commutative partially ordered group with two partial orderings. If the mixed envelopes $x\uenv y$ and $x\lenv y$ exist in $G$ for all $x,y\in G$ then $G$ is called a \emph{mixed lattice group}. Groups with mixed lattice order structure were first studied by Eriksson-Bique in \cite{eri1}, while the general definition of a mixed lattice group was introduced in \cite{eri}.

The following identities hold in any mixed lattice group. 

\begin{equation}\label{p1}
x\uenv y \, + \, y \lenv x \, = \, x+y
\end{equation}
\begin{equation}\label{p2}
z \,+ \,x\uenv y \, = \, (x+z)\uenv(y+z) \, \textrm{ and } \, z \,+\, x\lenv y \, = \, (x+z)\lenv(y+z)
\end{equation}
\begin{equation}\label{p3}
-(x\uenv y) \, = \, (-x) \lenv (-y)
\end{equation}

For proofs of these identities we refer to \cite{eri}.

The next lemma gives a useful criterion for checking that a group with two orderings is a mixed lattice group, and many of the examples that follow can easily be seen to be mixed lattice groups with this result. The proof can be found  in \cite[Theorem 3.3]{eri}. 

\begin{lemma}\label{mlg_lemma}
Let $G$ be a commutative partially ordered group with two partial orderings. Then $G$ is a mixed lattice group if and only if the element $x\uenv 0$ exists for all $x\in G$ (or equivalently, the element $0\uenv x$ exists for all $x\in G$).
\end{lemma}

We will use the following notation for the sets of positive elements. If $G$ is a mixed lattice group then we define $G_p=\{w\in G: w\geq 0\}$ and $G_{sp}=\{w\in G: w\sgeq 0\}$. If $S$ is a subset of $G$ then the sets of positive elements of $S$ are denoted by $S_p = S\cap G_p$ and $S_{sp}=S\cap G_{sp}$. 

We note that, unless $G=\{0\}$, $\leq$ and $\sleq$ cannot be trivial orderings (i.e. $G_{sp}\neq \{0\}$ and $G_{p}\neq \{0\}$) since the mixed envelope $x\uenv y$ must exist for all $x,y\in G$. Moreover, $G_{sp}\cap G_p \neq \{0\}$. This can be seen by taking nonzero elements $x\geq 0$ and $y\sgeq 0$ and noting that $0\leq x\leq y\uenv x$ and $0\sleq y\sleq y\uenv x$. Hence $y\uenv x$ is a nonzero element that belongs to $G_{sp}\cap G_p$.

Next we shall examine subsemigroups of a mixed lattice group. Our main interest is in those sub-semigroups that have a mixed order structure that was studied in the preceding section. We begin with a simple lemma.

\begin{lemma}\label{definitions2}
Let $G$ be a mixed lattice group and $S$ a sub-semigroup of $G_{sp}$ or $G_{p}$ such that $0\in S$.  The following statements are equivalent.
\begin{enumerate}[(i)]
\item
$M=(S,\leq, \sleq)$, where $\leq$ and $\sleq$ are the partial orderings of $G$ restricted to $S$,  
is a partially ordered semigroup in which the mixed envelopes exist and they are the same as in $G$. 
\item
If $S$ is considered as a subset of $G$, then the elements  
$x\lenv y$ and $x\uenv y$ belong to $S$ for all $x,y\in S$. 
\end{enumerate}
\end{lemma}

\begin{bproof}
It is clear that (i) implies (ii). If (ii) holds then we need to show that the mixed envelopes exist in $M$  and they are the same as in $G$. Let $\uenv$ and $\lenv$ be the mixed envelopes in $G$, and let the mixed envelopes in $M$ be denoted by $\uenv_{\!\!\! \footnotesize{M}}$ and $\lenv_{\!\! \footnotesize{M}}$  (whenever they exist). 
If $x,y\in S$ then by assumption, $w=x\uenv y \in S$ and $x\sleq w$ and $y\leq w$ hold in $M$. If $z\in S$ is any element such that $x\sleq z$ and $y\leq z$ in $M$ then these inequalities hold also in $G$. Hence, $w\leq z$ holds in $G$ and therefore it holds in $M$. This shows that $w = x \uenv_{\!\!\! M} \,y$, and so $\uenv_{\!\! M}$ exists and it is the same as $\uenv$. A similar reasoning shows that $\lenv_{\!\! M}$ exists too and it is the same as $\lenv$. Hence, (i) holds.
\end{bproof}

%---------------

The characterizations of Theorems \ref{karakt11} and \ref{karakt22} motivate the following definition. 

\begin{defn}\label{substructures}
Let $G$ be a mixed lattice group. 
\begin{enumerate}[(i)]
\item
A sub-semigroup $S$ of $G_{sp}$ containing $0$ is called an \emph{$A$-substructure} of $G$ if $x\uenv y - x\in S$ for all $x,y\in S$. 
\item
A sub-semigroup $S$ of $G_{p}$ containing $0$ is called a \emph{$B$-substructure} of $G$ if $x\uenv y - y\in S$ for all $x,y\in S$.
\item
An $A$-substructure $S$ of $G$ is called a \emph{mixed lattice sub-semigroup} of $G$ if 
$x\lenv y \in S$ for all $x,y\in S$. 
\end{enumerate}
\end{defn}

\begin{thm}\label{substructures2}
Let $G$ be a mixed lattice group and $S$ a sub-semigroup of $G$. Consider the partially ordered group 
$M=(S,\leq, \sleq)$, where $\leq$ and $\sleq$ are the partial orderings of $G$ restricted to $S$. 
Then the following hold. 
\begin{enumerate}[(a)]
\item
If $S$ is an $A$-substructure of $G$ then $M$ is an $A$-structure in its own right.
\item
If $S$ is a $B$-substructure of $G$ then $M$ is a $B$-structure in its own right.
\item
If $S$ is a mixed lattice sub-semigroup of $G$ then $M$ is a mixed lattice semigroup in its own right. 
\end{enumerate}
\end{thm}

\begin{bproof}
(a)\, If $S$ is an $A$-substructure of $G$ then $M$ is naturally ordered with respect to $\sleq$, for if $x,y\in S$ then  $x\sleq y$ is equivalent to $x\uenv y=y$ (by \eqref{p5b}). This is equivalent to $y-x=x\uenv y -x \in S$. Hence  $x\sleq y$ if and only if $y=x+w$ for some $w\in S$ (namely $w=y-x$), and thus $M$ is naturally ordered with respect to $\sleq$. Moreover, the cancellation property \eqref{canc1} holds in $M$ since it holds in $G$. Finally, since $x\uenv y -x \in S$ it follows that $x+(x\uenv y -x)=x\uenv y \in S$ for all $x,y\in S$. Thus the mixed upper envelopes exist in $M$ by Lemma \ref{definitions2}, and so $M$ is an $A$-structure. 
Part  (b) is proved by a similar argument.  
Part (c) then follows from (a) by noting that if $y\lenv x$ exists in $S$ then by \eqref{p1} the identity $x\uenv y + y\lenv x = x+y$ holds in $S$, and therefore it holds in $M$, by Lemma \ref{definitions2}.
\end{bproof}

The significance of $A$-structures in the theory of mixed lattice groups is largely due to the following result (\cite[Theorem 3.4]{eri}). 

\begin{thm}\label{astructure_in_mlgroup}
If $G$ is a mixed lattice group then the set $G_{sp}=\{w\in G: w\sgeq 0\}$ is an $A$-structure.  
On the other hand, a group extension of an $A$-structure is a mixed lattice group. 
\end{thm}

\begin{remark} 
By a group extension we mean the standard procedure of extending a semigroup to a group of differences by considering the equivalence classes of elements (c.f. \cite[Theorem II.3.5]{hebi}). For details of this construction applied to mixed lattice semigroups, see \cite{eri1}.
\end{remark}

The first part of Theorem \ref{astructure_in_mlgroup} follows also from the following characterization of mixed lattice groups. It basically states that $G_{sp}$ is an $A$-substructure and $G_p$ is a $B$-substructure.

\begin{thm}\label{mlg_char1}
Let $G=(G,+,\leq,\sleq)$ be a partially ordered commutative group with two partial orderings. The following conditions are equivalent.
\begin{enumerate}[(a)]
\item
$G$ is a mixed lattice group.
\item
The element $m=\min\{w\sgeq 0:y\leq w+x\}$  
exists for all $x,y\in G$, and it is given by $m=x\uenv y -x=y-y\lenv x$.  
\item
The element $r=\min\{w\geq 0:y\sleq w+x\}$ 
exists for all $x,y\in G$, and it is given by $r=y-x\lenv y = y\uenv x -x$
\end{enumerate}
\end{thm}

\begin{bproof} 
Let $x,y\in G$ and $A=\{w\sgeq 0:y\leq w+x\}$. For any $w\in A$ we set $z=w+x$, or $w=z-x$. Then %the element
$$
\begin{array}{rcl}
\min A &=& \min \{w: w\sgeq 0 \,\textrm{ and }\, y\leq w+x\} \\[1ex]
& = & \min \{w+x: w\sgeq 0 \,\textrm{ and }\, y\leq w+x\}-x \\[1ex]
& = & \min \{z: z-x\sgeq 0 \,\textrm{ and }\, y\leq z\}-x \\[1ex]
&=&  \min \{z: x\sleq z \,\textrm{ and }\, y\leq z\}-x \\[1ex]
&= & x\uenv y -x.
\end{array}
$$
Hence $\min A$ exists if and only if $x\uenv y -x$ exists for all $x,y\in G$, and the latter holds if and only if $G$ is a mixed lattice group. This also shows that $\min A=x\uenv y -x=y-y\lenv x$, where the last equality follows by \eqref{p1}. 
The equivalence (a)$\iff$(c) is proved similarly.
\end{bproof}

The equivalence (a)$\iff$(b) of the above theorem was essentially proved in \cite{eri} (Theorems 2.2 and 3.4) but the proofs and the statements of their results were different than in our Theorem \ref{mlg_char1}. 

Now if we consider the semigroup $G_{sp}$ with the partial orderings inherited from $G$ then, since $G_{sp}$ is the positive cone of $G$ for $\sleq$, it is naturally ordered with respect to $\sleq$. By part (b) of the preceding theorem, the element $\min\{w\in G_{sp} :y\leq w+x\}=x\uenv y-x$ exists in $G_{sp}$ for all $x,y\in G_{sp}$. By Theorem \ref{substructures2} this shows that $G_{sp}$ is an $A$-structure. Similarly, the set $G_p$ is naturally ordered with respect to $\leq$, and it follows from Theorem \ref{substructures2} and part (c) of Theorem \ref{mlg_char1} that $G_p$ is a $B$-structure. 
We have therefore the following:

\begin{cor}\label{positive_cones}
If $G$ is a mixed lattice group then $G_{sp}$ is an $A$-substructure and $G_p$ is a $B$-substructure.  
\end{cor}

\begin{remark}
The reason we haven't introduced a notion of a weak mixed lattice sub-semigroup of $G$ is that by Theorem \ref{wmlsg_is_positive}, a weak mixed lattice semigroup is a positive $A$-structure, so if $G$ contains a weak mixed lattice semigroup $S$ as a sub-semigroup of $G$ then necessarily $x\uenv y-x\in S$ and $x\lenv y\in S$ for all $x,y\in S$.   
By Definition \ref{substructures}, $S$ is then a mixed lattice sub-semigroup of $G$. Hence the notion of a weak mixed lattice sub-semigroup coincides with that of a mixed lattice sub-semigroup.
\end{remark}

The following examples clarify the notions introduced above and the importance of Lemma \ref{definitions2}.

\begin{example}\label{weak_mlsg_ex2}
Consider the additive group $\R$. Let $\leq$ be the usual order of $\R$, and define $\sleq$ as the partial ordering with the positive set $S=\{0\}\cup [2,\infty)$. Then $G=(\R,\leq,\sleq)$ is a mixed lattice group and $S$ is a sub-semigroup of $G$.  
Let us then consider the partially ordered semigroup $M=(S,\leq, \sleq)$ with  
the partial orderings $\leq$ and $\sleq$ inherited from $G$. Since $S=G_{sp}$, it follows by Corollary \ref{positive_cones} that $S$ is an $A$-substructure of $G$, and $M$ is then a positive $A$-structure, by Theorem \ref{substructures2}. It was also noted in \cite[Example 2.3]{jj1} that the mixed lower envelopes exist in $M$ for all $x,y\in M$. 
Thus $M$ is a weak mixed lattice semigroup. However, it was shown in \cite[Example 2.3]{jj1} that the identity $x\uenv y+y\lenv x=x+y$ does not hold in $M$. Therefore $M$ is not a mixed lattice semigroup, and $S$ is not a mixed lattice sub-semigroup of $G$. 
The crucial point is that the mixed envelopes in $M$ are not the same as the mixed envelopes in $G$. For instance, we have $3\lenv 2=0$ in $M$, but $3\lenv 2=1$ in $G$, and so $3\lenv 2 \notin S$.  
Note also that the mixed lower envelope in $M$ is not translation invariant. 
For example, $3\lenv 2 + 2 =2$, since $3\lenv 2=0$ in $M$, but $(3+2)\lenv (2+2)=5\lenv 4=3$. Hence $(3+2)\lenv (2+2)\neq 3\lenv 2 +2$. 
The subset $\{0,2,4,6,\ldots\}$ is a an example of a mixed lattice sub-semigroup of $G$. 
\end{example}

\begin{remark}\label{huom1} 
By Lemma \ref{definitions2} we might be tempted to conclude that if $S$ is a sub-semigroup of $G_{sp}$ such that $0\in S$ and $x\uenv y\in S$ whenever $x,y\in S$, then $S$ is an $A$-structure in its own right. However, this is not the case because such $S$ is not necessarily naturally ordered by the partial ordering $\sleq$ inherited from $G$.  
In fact, the condition $x\uenv y -x \in S$ for all $x,y\in S$ implies that $S$ is naturally ordered with respect to $\sleq$, as we showed in the proof of Theorem \ref{substructures2}.   
Observe also that $x\uenv y -x \in S$ implies that $x\uenv y\in S$, but not conversely. Similar remarks apply for  
$B$-substructures. An example of this is provided next.
\end{remark}

\begin{example}\label{huom1ex}
Let $G=(\Z\times \Z,+,\leq,\sleq)$, where $+$ is the usual coordinatewise addition, and $\leq$ and $\sleq$ are the  lexicographic orderings defined by
$$
x\leq y \quad \iff \quad  ( \; x_1< y_1 \quad \textrm{or} \quad x_1=y_1 \quad \textrm{and} \quad x_2\leq y_2 \;)
$$
and
$$
x\sleq y \quad \iff \quad ( \;x_2< y_2 \quad \textrm{or} \quad x_2=y_2 \quad \textrm{and} \quad x_1\leq y_1 \;).
$$
Then $G$ is a mixed lattice group (c.f. \cite[Example 2.25]{jj1}). 
Now the sub-semigroup $F=G_{sp} \cap G_p$ has the property that $x\uenv y\in F$ for all $x,y\in F$ but $F$ is not naturally ordered, and hence the element $x\uenv y -x$ is not in $F$ for all $x,y\in F$. For example, choose $x=(0,1)$ and $y=(1,1)$. Then $y\uenv x =(0,2)\in F$ but $y\uenv x-y=(-1,1)\notin F$.   
\end{example}

A subgroup $E$ of $G$ is called a \emph{mixed lattice subgroup} of $G$ if $x\uenv y \in E$ and $x\lenv y\in E$ whenever $x,y\in E$. 
The ordered sub-semigroup structures discussed above are closely related to the structure of mixed lattice subgroups. 

\begin{thm}\label{mlsubgroup1}
If $E$ is a mixed lattice subgroup of $G$ then $E_{sp}$ is an $A$-substructure of $G$ and $E_{p}$ is a $B$-substructure of $G$. Conversely, if $S$ is an $A$-substructure or a $B$-substructure of $G$ then the subgroup generated by $S$ is a mixed lattice subgroup of $G$.
\end{thm}

\begin{bproof}
If $E$ is a mixed lattice subgroup then $E_{sp}$ is clearly a sub-semigroup of $G_{sp}$ and $0\in E_{sp}$. If $x,y\in E_{sp}$ then $x\uenv y\in E$ and $0\sleq x\uenv y -x \in E$. Hence $x\uenv y -x \in E_{sp}$, proving that $E_{sp}$ is an $A$-substructure of $G$. A similar argument shows that $E_{p}$ is a $B$-substructure of $G$.

For the second part, let $S$ be an $A$-substructure of $G$ and let $E=S-S$ be the subgroup generated by $S$. If $x,y\in E$ then $x=u-v$ and $y=a-b$ for some $u,v,a,b\in S$. Then by the identity \eqref{p2} we obtain 
$$
x\uenv y = (u-v)\uenv(a-b) = (u+b)\uenv (a+v) - (v+b) \in S-S.
$$ 
This implies that $E$ is a mixed lattice subspace. 
A similar argument applies in the case that $S$ is a $B$-substructure.
\end{bproof}

Mixed lattice groups have the following classification according to the structure of the sets of positive elements. 

\begin{defn}\label{mlgroups}
Let $G$ be a mixed lattice group. Then $G$ is called 
\begin{enumerate}[(i)]
\item
\emph{pre-normal} if $(-G_{sp})\cap G_{p}=\{0\}$.
\item
\emph{normal} if $G_{sp}\subseteq G_p$ or $G_p\subseteq G_{sp}$.
\item
\emph{pre-regular} if $G_{sp}\subseteq G_p$.   
\item
\emph{quasi-regular} if $G_{sp}$ is a mixed lattice sub-semigroup of $G$.
\item
\emph{regular} if $G$ is quasi-regular and $G_{sp}$ is generating (i.e. $G=G_{sp}-G_{sp}$).
\end{enumerate}
\end{defn}

The notions of regular, quasi-regular and pre-regular mixed lattice groups were introduced in \cite{jj1}.  
It was proved in \cite[Lemma 2.4]{jj5} that 
$G$ is pre-normal if and only if the following equivalences hold in $G$:

\begin{equation}\label{r0}
x\leq y \, \iff \, y\uenv x = y \, \iff \, x\lenv y = x. \\[2ex]
\end{equation}

Moreover, it is clear that if $G_{sp}\subseteq G_p$ or $G_p\subseteq G_{sp}$ holds then $(-G_{sp})\cap G_{p}=\{0\}$. 
Hence every normal mixed lattice group is pre-normal. 
By Theorem \ref{astructure_in_mlgroup}, 
pre-regularity is equivalent to the condition that $G_{sp}$ is a positive $A$-structure, and it clearly implies normality. Quasi-regularity in turn implies pre-regularity, by Theorem \ref{wmlsg_is_positive}. 

None of these implications hold conversely (see also examples in \cite{jj1}).  

\begin{example}
To see that pre-normality does not imply normality, let $G=(\R^3,+,\leq,\sleq)$ with the usual addition and the positive cones $G_p=\{a(1,0,1)+b(-1,0,1):a,b\geq 0\}$ and $G_{sp}=\{a(0,1,1)+b(0,-1,1):a,b\geq 0\}$. Then $G$ is a pre-normal mixed lattice group but not normal. This example shows also that generally neither of the positive cones $G_p$ or $G_{sp}$ need to be generating.
\end{example}

\begin{example}\label{uusi_irregular}
Let $G=(\R^2, \leq, \sleq)$, where $\leq$ is the lexicographic ordering defined as in Example \ref{huom1ex}, and $\sleq$ is the partial order with the positive cone $G_{sp}=\{\alpha (-1,1)+ \beta (1,1): \alpha,\beta \geq 0\}$. Then $G$ is a mixed lattice group that is not pre-normal.  
For example, $x=(-1,1)\in G_{sp}\cap -G_p$, and hence $G_{sp}\cap -G_p \neq \{0\}$. 
Here $G_{sp}$ is another example of an $A$-structure that is not a positive $A$-structure. 
For example, if $u=(-1,1)$ and $v=(1,1)$ then $u\uenv v = (1,3)$ and $u+v=(0,2)$, so $v\uenv u \leq u+v$ does not hold (c.f. Theorem \ref{positive_astr_char}). 
\end{example}

Quasi-regular mixed lattice groups can be characterized in many ways, some of which are collected in the following theorem. For other equivalent conditions, see \cite[Theorem 2.7]{jj5}.

\begin{thm}\label{q-regular_char}
Let $G$ be a mixed lattice group and $x,y,z,u\in G$. The following statements are equivalent.
\begin{enumerate}[(a)]
\item
$G$ is quasi-regular.
\item
If \,$u\sleq x\sleq z$ \,and \, $u\sleq y\sleq z$\, then\, $x\uenv y \sleq z$\, and\, $u \sleq x\lenv y$.
\item
If \,$x\sleq y$ \,then\, $z\lenv x \sleq z\lenv y$\, and\, $z\uenv x \sleq z\uenv y$ \,for all $z\in G$.
\item
The distributive inequalities \,$(z\lenv x)\uenv (z\lenv y) \leq z\lenv(x\uenv y)$ \,and \, $z\uenv(x\lenv y) \leq (z\uenv x)\lenv (z\uenv y)$ 
hold in $G$. 
\item
The associative inequalities \, $x\lenv (y\lenv z) \leq (x\lenv y)\lenv z$ \, and \\ 
$(x\uenv y)\uenv z \leq x\uenv (y\uenv z)$ \,
hold in $G$. 
\item
The modular inequalities $(z\lenv x)\uenv (z\lenv y) \leq z\lenv((z\lenv x)\uenv y)$ \, and \, $(z\uenv x)\lenv (z\uenv y) \geq z\uenv((z\uenv x)\lenv y)$ 
hold in $G$.
\end{enumerate}
\end{thm}

\begin{bproof}
The equivalence (a)$\iff$(b) was proved in \cite[Theorem 3.5]{eri}, the equivalence (b)$\iff$(f) in \cite[Theorem 2.7]{jj5}, the implication (a)$\implies$(c) in \cite[Theorem 3.2]{jj2}, and  
(c)$\implies$(b) in \cite[Theorem 2.8]{jj5}. It follows from the identity \eqref{p3} that the two distributive inequalities are equivalent, and similarly for the equivalence of the two associative inequalities. Then the equivalences (c)$\iff$(d)$\iff$(e) follow from \cite[Theorem 2.8]{jj5}.
\end{bproof}

For mixed lattice subgroups, the next result follows easily from Theorems \ref{mlsubgroup1} and \ref{q-regular_char}(b). 

\begin{thm}
If $G$ is quasi-regular and $E$ is a mixed lattice subgroup of $G$ then $E_{sp}$ is a mixed lattice sub-semigroup of $G$.
\end{thm}

\begin{bproof}
By Theorem \ref{mlsubgroup1}, $E_{sp}$ is an $A$-substructure, and if $x,y\in E_{sp}$ then, since $E$ is a mixed lattice subgroup, we have $x\lenv y \in E$. But from $x\sgeq 0$ and $y\sgeq 0$ it follows by Theorem \ref{q-regular_char}(b) that $x\lenv y \sgeq 0$, so that $x\lenv y \in E_{sp}$. Hence $E_{sp}$ is a mixed lattice sub-semigroup of $G$.
\end{bproof}

The notion of regularity can be extended to include single elements.

\begin{defn}
An element $x$ of a mixed lattice group $G$ is called a \emph{regular element} if $x$ can be written as $x=u-v$ with $u,v\in G_{sp} \cap G_p$.
\end{defn}

With the above definition, a quasi-regular mixed lattice group $G$ is regular if and only if every $x\in G$ is a regular element. The notion of a regular element will be useful in the next section.

\begin{thm}
The set of regular elements is a mixed lattice subgroup.   
\end{thm}

\begin{bproof}
Let $F=G_{sp} \cap G_p$. Then $F$ is clearly a sub-semigroup of $G$ and the set $E=F-F$ of all regular elements is a subgroup. We first note that if $u,v\in F$ then $0\sleq u\sleq u\uenv v$ and $0\leq v\leq u\uenv v$, so $u\uenv v\in F$. Now if $x,y\in E$ then $x=u-v$ and $y=a-b$ for some $u,v,a,b\in F$. Then, as in the second part of the proof of Theorem \ref{mlsubgroup1}, we conclude that 
$$
x\uenv y = (u-v)\uenv(a-b) = (u+b)\uenv (a+v) - (v+b) \in F-F.
$$ 
Hence $E$ is a mixed lattice subgroup. 
\end{bproof}

\begin{remark}
We should point out that the set $F=G_{sp} \cap G_p$ is not necessarily an $A$-substructure of $G$ (see Remark \ref{huom1} and Example \ref{huom1ex}).
\end{remark}

\section{Archimedean mixed lattice groups}\label{sec:archi}

In this section we study the implications of the Archimedean property in mixed lattice groups. 
A mixed lattice  group $G$ is called $(\leq)$-\emph{Archimedean} if the condition $nx\leq y$ for all $n\in \N$ implies $x\leq 0$. Similarly, $G$ is called $(\sleq)$-\emph{Archimedean} if the condition $nx\sleq y$ for all $n\in \N$ implies $x\sleq 0$. It is easy to see that if $G$ is pre-regular and 
$(\leq)$-Archimedean then it is also $(\sleq)$-Archimedean, 
but the converse is not true (see \cite[Example 5.20]{jj2}). Similarly, if $G_p \subseteq G_{sp}$ then the $(\sleq)$-Archimedean property implies that $G$ is $(\leq)$-Archimedean.

We begin with the following result on regular $(\leq)$-Archimedean mixed lattice groups. %The non-trivial parts
The next theorem was partially proved in \cite[Theorem 4.1]{eri1}.

\begin{thm}\label{archim_theorem}
In a mixed lattice group $G$ the following properties are equivalent.
\begin{enumerate}[(a)]
\item
If $x\lenv y=0$ then $(mx)\lenv (ny)=0$ for all $m,n\in\N$.
\item
The identity $(nx)\lenv (ny)=n(x\lenv y)$ holds in $G$ for all $n\in \N$. %
\item
The identity $(nx)\uenv (ny)=n(x\uenv y)$ holds in $G$ for all $n\in \N$. %
\end{enumerate}
Moreover, the properties (a), (b) and (c) hold in every regular $(\leq)$-Archimedean mixed lattice group.  
\end{thm}

\begin{bproof}
The implication (a)$\implies$(b) was established in the proof of \cite[Theorem 4.1]{eri1} for $m=1$. 
Assume that (b) holds and $x\lenv y=0$. Then $x\sgeq 0$ and $y\geq 0$, and consequently, $mx\sgeq 0$ and $ny\geq 0$ for all $m,n\in \N$. Then $0\leq (mx)\lenv (ny) \leq (px)\lenv (py) = p(x\lenv y) =0$, where $p=\max\{m,n\}$. Hence  
$(mx)\lenv (ny)=0$ for all $m,n\in\N$. 

If (b) holds then by \eqref{p1} we have
$$
(nx)\uenv (ny)=nx +ny - (ny)\lenv (nx)=n(x+y-y\lenv x)=n(x\uenv y),
$$
and the implication (c)$\implies$(b) is proved in the same manner.  
The fact that (a), (b) and (c) hold in regular $(\leq)$-Archimedean mixed lattice groups was proved in \cite[Theorem 4.1]{eri1}. 
\end{bproof}

It follows from part (b) of Theorem \ref{archim_theorem} that if $G$ is a regular $(\leq)$-Archimedean mixed lattice group then $nx\leq ny$ implies $x\leq y$ (this was also proved in \cite[Theorem 4.1]{eri1}). 
This also holds in every commutative lattice ordered group (\cite[XIII.3, Corollary 2]{birk}), with or without the Archimedean property.  
However, we cannot drop the Archimedean and regularity assumptions without losing this property as the following example shows.

\begin{example} %
Let $G=(\Z \times \Z,+,\leq,\sleq)$ where $+$ is the usual coordinatewise addition. We define the partial orderings $\leq$ and $\sleq$ on $G$ by 
$$
(x,y)\geq (0,0) \;\, \iff \;\, (x,y)=(0,0) \;\, \text{or} \;\, (x,y)=(1,1) \;\, \text{or} \;\,(x\geq 2 \; \text{and}\; y\geq 2)
$$
and
$$
(x,y)\sgeq (0,0) \quad \iff \quad x\geq 0, \, y\geq 0 \;\, \text{ and } \;\, x=y.
$$
Then $G$ is a quasi-regular mixed lattice group but not $(\leq)$-Archimedean (for example, if $x=(0,-1)$ and $y=(2,2)$ then $nx\leq y$ for all $n\in \N$ but $x\leq 0$ does not hold).  
Now $nu\geq nv$ does not imply that $u\geq v$. This can be seen by choosing $u=(1,2)$. Then $2u\geq 0$ but $u\geq 0$ does not hold. 
\end{example}

Similarly, $nx\sgeq ny$ does not necessarily imply $x\sgeq y$  
even if $G$ is a $(\sleq)$-Archimedean quasi-regular mixed lattice group (see Example \ref{uusiesim}).

It is also known that a lattice ordered group cannot contain nonzero elements of finite order (\cite[XIII.3, Corollary 1]{birk}).  
In mixed lattice groups the situation is more complicated but 
we will now show that nonzero elements of finite order cannot exist in a rather general class of normal mixed lattice groups with the Archimedean property.

\begin{thm}\label{inforder}
Let $G$ be a normal mixed lattice group such that one of the following conditions hold.
\begin{enumerate}[(a)]
\item
$G_{sp}\subseteq G_p$ and $G$ is $(\leq)$-Archimedean.
\item
$G_{p}\subseteq G_{sp}$ and $G$ is $(\sleq)$-Archimedean.
\end{enumerate}
Then every nonzero element of $G$ has infinite order.
\end{thm}

\begin{bproof}
Assume that an element $x\neq 0$ has order $k>1$. Then $kx=0$ and every $n\in \N$ can be written as $n=mk+p$ for some $m\in\{0,1,2,\ldots\}$ and $p\in\{0,1,\ldots, k-1\}$. Thus $nx=mkx+px=px$. 
Assume first that the condition (a) holds. 
If $p=0$ then $nx=0$, and if $p\neq 0$ then, since $x\uenv 0 \geq 0$ we have 
$$
nx=px\sleq p(x\uenv 0) \leq (k-1)(x\uenv 0). 
$$
By (a), this implies that $nx\leq (k-1)(x\uenv 0)$ for all $n$, and hence $x\leq 0$. But then 
$0=kx\leq (k-1)x\leq \ldots \leq x\leq 0$, and this implies that $x=0$.

Suppose next that the condition (b) holds. 
If $p=0$ then $nx=0$, and if $p\neq 0$ then, since $0\uenv x \sgeq 0$ we have 
$$
nx=px\leq p(0\uenv x) \sleq (k-1)(0\uenv x). 
$$
By (b), this implies that $nx\sleq (k-1)(0\uenv x)$ for all $n$, and hence $x\sleq 0$. But then 
$0=kx\sleq (k-1)x\sleq \ldots \sleq x\sleq 0$, and this implies that $x=0$. 

Hence, in both cases $kx=0$ for some $k\in \N$ implies that $x=0$. Thus every nonzero element has infinite order. 
\end{bproof}

On the other hand, we have the following result regarding elements of finite order.

\begin{thm}\label{regular_inforder}
Let $G$ be a mixed lattice group. If $G$ is $(\leq)$-Archimedean or $(\sleq)$-Archimedean and $x\in G$ is a nonzero element of finite order then $x$ cannot be a regular element.
\end{thm}

\begin{bproof}
Let $x\neq 0$ be an element of order $k$. Then, as shown in the proof of Theorem \ref{inforder}, for every $n\in\N$ we have $nx=px$ for some $p\in\{0,1,\ldots,k-1\}$. 
If $x$ is a regular element, then there exists $u\in G_{sp} \cap G_p$ such that $x\sleq u$ and $x\leq u$. 
Assume that $G$ is $(\sleq)$-Archimedean. Then for every $n\in \N$ there is some $p\in\{0,1,\ldots,k-1\}$ such that 
$nx=px\sleq pu \sleq (k-1)u$. By the Archimedean property, this implies that $x\sleq 0$, and consequently, $x=0$. This is a contradiction, and so $x$ cannot be a regular element. 

Since we also have $x\leq u$, we can use a similar reasoning to prove the result in the $(\leq)$-Archimedean case.  
\end{bproof}

Combining the results of the preceding theorems yields the following:

\begin{cor}\label{inforder_cor}
A mixed lattice group $G$ does not contain nonzero elements of finite order if one of the following conditions hold.
\begin{enumerate}[(a)]
\item
$G_{sp}\subseteq G_p$ and $G$ is $(\leq)$-Archimedean.
\item
$G_{p}\subseteq G_{sp}$ and $G$ is $(\sleq)$-Archimedean.
\item
$G$ is a $(\leq)$-Archimedean or $(\sleq)$-Archimedean mixed lattice group such that the set $G_{sp}\cap G_p$ is generating.
\end{enumerate}
\end{cor}

The preceding result does not hold without the Archimedean hypothesis even under the stronger quasi-regularity assumption. In the following example we construct a quasi-regular non-Archimedean mixed lattice group that contains elements of finite order. This also provides a counterexample for some other remarks we have made earlier.

\begin{example}\label{uusiesim}
Let $F=(\Z,+)$ with the usual addition and ordering of integers, and let $H=(\{0,1\},\oplus)$ where the addition is defined by the rules $0\oplus 0=0=1 \oplus 1$ and $0\oplus 1=1\oplus 0=1$. If $H$ is equipped with the trivial ordering (i.e. $x\leq y$ iff $x=y$)   
then $H$ is a partially ordered group. Let $X=F\times H$ and define addition $\boxplus$ on $X$ by 
$$
(x_1 ,y_1) \boxplus (x_2 , y_2)= (x_1+y_1, x_2 \oplus y_2)
$$ 
and let the partial order $\leq$ be defined by
$$
(0,0)\leq (x,y)  \iff 0<x \quad \text{or} \quad (x,y)=(0,0). %
$$
Then $(X,\boxplus,\leq)$ is a partially ordered group. Next, we define another partial ordering $\sleq$ on $X$ by 
$$
(x,y)\sgeq (0,0) \iff (\, y=0 \;\, \text{and} \;\, x=2k \,) \quad \text{or} \quad
(\, y=1 \;\, \text{and} \;\, x=2k+1 \,),
$$
where $k=0,1,2,\ldots$. 
Then $G=(X,\boxplus,\leq,\sleq)$ is a mixed lattice group where the element $(0,1)$ has order 2. 
That $G$ is a mixed lattice group is easy to check using Lemma \ref{mlg_lemma},  
and checking that $x\uenv 0$ exists in all different cases $x\leq 0$, $x\sgeq 0$, $x=(2k,1)$ and $x=(2k+1,0)$. 
Moreover,  we have 
$$
(0,0)\sleq (1,1)\sleq(2,0) \sleq (3,1) \sleq \ldots \sleq (2k,0) \sleq (2k+1,1) \sleq  \ldots,
$$
and the same chain of inequalities hold also for the ordering $\leq$. Hence the set $G_{sp}$ is totally ordered with respect to both orderings. This implies that if $x,y\in G_{sp}$ then $x\lenv y =\min \{x,y\}\in G_{sp}$, and so $G$ is quasi-regular. 

We also observe that $G$ is $(\sleq)$-Archimedean but not $(\leq)$-Archimedean. We first note from the above chain of inequalities that $(2k,0)=2k(1,1)$ and $(2k+1,1)=(2k+1)(1,1)$ hold for all $k\in \Z$. Hence $x\sleq 0$ or $x\sgeq 0$ if and only if $x$ is a multiple of the element $(1,1)$. In this case, 
$y\sleq x$ implies that $y=p(1,1)$ for some $p\in \Z$. It follows that if $ny\sleq x$ for all $n\in\N$ then necessarily $y\sleq 0$. On the other hand, if neither $x\sleq 0$ nor $x\sgeq 0$ holds then $x=(2k,1)$ or $x=(2k+1,0)$ for some $k\in \Z$. In the first case, if $y\sleq x$ then $y=(2p,1)$ or $y=(2p-1,0)$ for some $p\leq k$. Now if $n$ is even, then in both cases we have $ny=(2m,0)$ for some $m\in \Z$. But then $ny\sleq x$ does not hold, and it is therefore not possible that $ny\sleq x$ for all $n\in \N$. The same is true if $x=(2k+1,0)$, by a similar reasoning. This shows that $G$ is $(\sleq)$-Archimedean. 
To see that $G$ is not $(\leq)$-Archimedean, we choose $x=(0,1)$. Then $nx\leq (1,0)$ for all $n\in\N$ but $x\leq 0$ does not hold. 

We already showed that $G$ is quasi-regular, but it is easy to see that $G$ is not regular (i.e. the set $G_{sp}\cap G_p$ is not generating). For example, if $x=(0,1)$ then there is no such element $u$ in $G$ such that $u\sgeq 0$ and $x\sleq u$. 
This example shows that both the the Arhimedean hypothesis and the assumption that the set $G_{sp}\cap G_p$ is generating are essential in Corollary \ref{inforder_cor} (a) and (c).

We shall make a couple of additional remarks regarding this example.  
First, we note that $nx\sgeq 0$ does not imply that $x\sgeq 0$. For example, $2(1,0)\sgeq (0,0)$ but $(1,0)\sgeq (0,0)$ does not hold.  
Secondly, $G$ is not a lattice with respect to either partial ordering. For example, the elements $\sup\{(0,0),(0,1)\}$ and $\spsup\{(0,0),(0,1)\}$ do not exist in $G$.  
\end{example}

%%-----------------------------
%seur

The theorems presented above can be used to improve certain results given in earlier papers. 
Before giving an example, we recall the following notation that was introduced in \cite{jj1}. The elements $\ul{x}=x\uenv 0 + 0\uenv (-x)$ and $\lu{x}=0\uenv x + (-x)\uenv 0$ are called the \emph{asymmetric absolute values} of $x$, as they possess similar basic properties as the usual absolute value in lattice ordered groups and Riesz spaces 
(for more on these, we refer to \cite{birk, lux}).

As an application, we consider mixed lattice group homomorphisms which were studied in \cite{jj1}. 
If $G$ and $H$ are mixed lattice groups then an additive mapping $T:G\to H$ is called a \emph{mixed lattice homomorphism} if $T(x\uenv y)=(Tx)\uenv(Ty)$ and $T(x\lenv y)=(Tx)\lenv(Ty)$ hold for all $x,y\in G$. 
It was shown in \cite[Corollary 5.6]{jj1} that the characterization of the following corollary holds if $H$ is a regular $(\leq)$-Archimedean mixed lattice group. With Corollary \ref{inforder_cor} we can state this result in a more general form as follows.

\begin{cor}\label{latticehomom3}
Let $T:G\to H$ be an additive mapping between two mixed lattice groups $G$ and $H$. If 
one of the conditions of Corollary \ref{inforder_cor} hold in $H$  
then the following statements are equivalent.
\begin{enumerate}[(a)]
\item \quad
$T$ is a mixed lattice homomorphism.
\item \quad
$T(\lu{x})=\lu{(Tx)}$ \;  for all \, $x\in G$.
\item \quad
$T(\ul{x})=\ul{(Tx)}$ \; for all \, $x\in G$.
\end{enumerate}
\end{cor}

\begin{bproof}
It was proved in \cite[Proposition 5.3]{jj1} that the given statements are equivalent if $H$ contains no nonzero elements of order 2. Thus the result follows from Corollary \ref{inforder_cor}.
\end{bproof}

%%------------------------------------------------------------------------------------------------

\bibliographystyle{plain}
%\bibliography{reference}

\end{document}